\documentclass[11pt,a4paper,reqno]{amsart}
\usepackage{graphicx}
\usepackage{epsfig}
\usepackage{amsmath}
\usepackage{amsfonts}
\usepackage{amssymb}
\usepackage{amscd}
\newtheorem{theorem}{Theorem}[section]
\newtheorem{corollary}[theorem]{Corollary}

\newtheorem{lemma}[theorem]{Lemma}
\newtheorem{proposition}[theorem]{Proposition}

\begin{document}

\begin{minipage}[b]{0.5\linewidth}
  {\includegraphics[height=0.72in,width=5.94in]{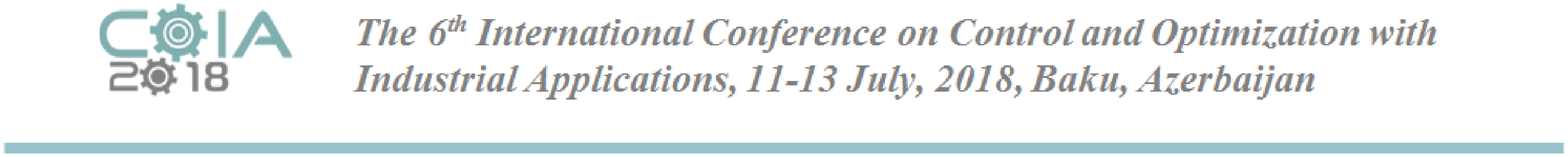}}
\end{minipage}

\bigskip
\title [R.P. Agaev, P.Y. Chebotarev: Models of latent consensus]
{TWO MODELS OF LATENT CONSENSUS IN MULTI-AGENT SYSTEMS}

\author[R.P. Agaev, P.Y. Chebotarev] {R. P. Agaev$^1$,  P. Y. Chebotarev$^{1,2}$}

\maketitle
\begin{center} $^{1\,}$Trapeznikov Institute of Control Sciences of RAS, 65 Profsoyuznaya Str., Moscow, Russia
\\ \indent $^2$Kotelnikov Institute of Radioengineering and Electronics of RAS, Moscow, Russia
\\e-mail: agaraf3@gmail.com, pavel4e@gmail.com \end{center}

\bigskip\bigskip
\section{Introduction} In this paper, we propose several consensus protocols of the first and second order for networked multi-agent systems and provide explicit representations for their asymptotic states. These representations involve the eigenprojection of the Laplacian matrix of the dependency digraph.
In particular, we study regularization models for the problem of coordination when the dependency digraph does not contain a converging tree. In such models of the first kind, the system is supplemented by a dummy agent, a ``hub" that uniformly, but very weakly influences the agents and, in turn, depends on them. In the models of the second kind, we assume the presence of very weak background links between the agents. Besides that, we present a description of the asymptotics of the classical second-order consensus protocol.

\medskip
\section{Notation and auxiliary results}
Consider the continuous protocol of consensus seeking in a multi-agent system \cite{Olfati-SaberMurray04,MesbahiEgerstedt10book}:
\begin{gather}\label{e_model0715}
\dot{x}_i(t)=-\sum_{j=1}^na_{ij}\left({x_i(t)-x_j(t)}\right),\quad
i=1,\ldots,n ,
\end{gather}
where $x_i(t)$ is the characteristic to be adjusted of agent~$i.$

With the \emph{dependency matrix\/} $A=(a_{ij})_{n\times n}$ of protocol \eqref{e_model0715} we associate the \emph{weighted dependency digraph} $\Gamma $ with vertex set $\{1,\ldots,n\}$ that for every element $a_{ij}>0$ of $A$ has an arc $(i,j)$ with weight $w_{ij}=a_{ij}$. Thus, arcs in $\Gamma $ are drawn in the direction of dependency. The weight $w_{ij}$ of arc $(i,j)$ reflects the strength of this influence.
The \emph{Laplacian matrix\/} of the weighted digraph $\Gamma $ is defined \cite{CheAga02ap} by
$$ 
L={\rm diag}(A  \textbf{1})-A.
$$ 
%
This matrix 
appears in the matrix representation of the model~\eqref{e_model0715}:
\begin{gather}\label{e_040715q}
\dot {\boldsymbol{x}}(t)=-L\,\boldsymbol{x}(t).
\end{gather}

The {\it eigenprojection\/} of matrix $A,$ {\it corresponding to the eigenvalue}~$0$ (or simply the {\it eigenprojection of\/}~$A$) is a projection (i.e., an idempotent matrix) $A^\vdash$ such that
${\mathcal{R}}(A^\vdash)={\mathcal{N}}(A^\nu)$ and
${\mathcal{N}}(A^\vdash)={\mathcal{R}}(A^\nu).$ Such a matrix $A^\vdash$ is always unique.

\begin{proposition}\label{p_limInf}
If $\boldsymbol{x}(t)$ is a solution to the system of equations \eqref{e_040715q}$,$ then
$$
\lim_{t\to\infty}\boldsymbol{x}(t)=L^\vdash  \boldsymbol{x}(0).
$$
\end{proposition}

\medskip
\section{Consensus protocols with a symmetric hub}
Consider the following protocols of consensus seeking:
\begin{gather}\label{e_modeLK040715c} \dot{
\boldsymbol{y}}(t)=
-(L_0+H_{\delta,\boldsymbol{v}})\,\boldsymbol{y}(t),
\end{gather}
where $\boldsymbol{y}(t)=(y_1(t),\ldots, y_{n+1}(t))^{\mathrm T}$, $H_{\delta  ,\boldsymbol{v}}=\delta   H_I+H_{\boldsymbol{v}}$,
$\delta  >0$,
$$ 
H_I =
\begin{pmatrix}
I       &- \textbf{1}\\[.3em]
\textbf{0}^{\mathrm T} &    0   \end{pmatrix}\quad \mbox{and}\quad
H_{\boldsymbol{v}} =
\begin{pmatrix}
0_{n\times n} & {\bf 1}\\[.3em]
   -\boldsymbol{v}^{\mathrm T}     &   s \end{pmatrix}
$$ 
are Laplacian matrices of order $n+1,\,$ $\boldsymbol{v}=(v_1,\ldots, v_n)^{\mathrm T},\,$ and $\,s=\sum_{i=1}^n v_i.$

Every protocol \eqref{e_modeLK040715c} is obtained from \eqref{e_040715q} by adding an ${(n+1)}$st agent that influences the $n$ agents of the system with intensity~$\delta  $ and is subjected to their influences whose intensities are given by the vector~$\boldsymbol{v}.$ This supplementary agent (and the corresponding vertex of the dependency digraph) will be called a \emph{hub}.
The motivation for the introduction of a hub is to assure the presence of a spanning in-tree in the dependency digraph (which guarantees asymptotic consensus) due to hub's uniform weak influences on the agents.

The eigenprojection $(L_0+H_{\delta  ,\boldsymbol{v}})^\vdash$ has a simple representation, which provides an expression for the consensus by means of Proposition~\ref{p_limInf}.

\begin{lemma}\label{02072015b}
$\,(L_0+H_{\delta  ,\boldsymbol{v}})^\vdash
={\bf 1}'\tfrac1{s+\delta  } \Big[\boldsymbol{v}^{\mathrm T}\big(I+\tfrac1\delta   L\big)^{-1}\;\;\, \delta  \Big].$
\end{lemma}

\begin{theorem}\label{t_UniHub}
In the system \eqref{e_modeLK040715c} with a hub that uniformly influences the agents$,$ the consensus is $\:\lim_{t\to\infty}\boldsymbol{y}(t)
={\bf 1}'\tfrac1{s+\delta  } \Big[\boldsymbol{v}^{\mathrm T}\big(I+\tfrac1\delta   L\big)^{-1}\;\; \delta \Big]\boldsymbol{y}(0).$
\end{theorem}

\section{Protocols with weak background links between the agents}
\label{s_backgr}

Now we turn to a different method of regularization. Consider the protocols of consensus seeking in which the initial dependency digraph $\Gamma $ is supplemented by the complete digraph of ``background links'' with uniform arc weights:
\begin{gather}\label{e_modeLK011115c}
\dot{ \boldsymbol{x}}(t)
=-(L+\delta K)\,\boldsymbol{x}(t),
\end{gather}
where $\delta>0$ and $K=I-E$ ($K$ is the Laplacian matrix of the complete digraph with arc weights~$1/n$).

First, in the same way as before, we study a more general type of regularization.
Let $\boldsymbol{v}$ be a distribution vector on the set of agents, i.e., ${v_i\ge 0}$ (${i=1,\ldots, n}$) and
${\sum_{i=1}^nv_i=1}$. Set
\begin{gather}\label{e_modeL-PR}
\dot{ \boldsymbol{x}}(t)=-(L+\delta   D)\,\boldsymbol{x}(t),
\end{gather}
where $\,\delta>0,$ $\,D=I-V,$ and $\,V={\bf 1}\boldsymbol{v}^{\mathrm T}.$ $D$ is Laplacian; $\delta
v_i$ is the influence of agent $i$ on any other agent.

\begin{lemma}\label{l_PRsop}
In the above notation\/$,$

$1.$ $ (L+\delta   D)^\vdash = {\bf 1} \boldsymbol{v}^{\mathrm T}
\big(I+\tfrac1\delta   L\big)^{-1}$$;$

$2.$ $ \lim_{\delta  \to0}(L+\delta   D)^\vdash ={\bf 1}
\boldsymbol{v}^{\mathrm T}\!\bar{J}.$
\end{lemma}

Proposition~\ref{p_limInf} and Lemma~\ref{l_PRsop} imply Theorem~\ref{t_limInf_0111a}.

\begin{theorem}\label{t_limInf_0111a}
If $\boldsymbol{x}(t)$ is a solution to the system of equations \eqref{e_modeL-PR}$,$ then 

$1.$ $\lim_{t\to\infty}\boldsymbol{x}(t) \!=\! {\bf 1}
\boldsymbol{v}^{\mathrm T} \big(I+\tfrac1\delta
L\big)^{-1}\boldsymbol{x}(0)$$;$
$2.$ $\lim_{\delta  \to+0}\lim_{t\to\infty}\boldsymbol{x}(t) ={\bf 1}
\boldsymbol{v}^{\mathrm T}\!\bar{J} \boldsymbol{x}(0).$
\end{theorem}

Observe that with ${\boldsymbol{v}=\tfrac1n {\bf 1}}$ we have ${{\bf 1}\boldsymbol{v}^{\mathrm T}=E}$ and ${D=K}$, and so Corollary~\ref{c_limInf} holds true.
\begin{corollary}\label{c_limInf}
If $\boldsymbol{x}(t)$ is the solution to the system of equations
\eqref{e_modeLK011115c}$,$ then

$1.$ $\lim_{t\to\infty}\boldsymbol{x}(t)
\!=\!{\bf 1}\,\tfrac1n \sum_{j=1}^ns_j^\delta    x_j(0),$ where
$s_j^\delta  $ are the column sums of $(I+\tfrac1\delta L\big)^{-1}$$;$

$2.$ $\lim_{\delta  \to+0}\lim_{t\to\infty}\boldsymbol{x}(t)
={\bf 1}\,\tfrac1n \sum_{j=1}^n\bar{J}_{ \cdot j}x_j(0),$ where
$\bar{J}_{\cdot j}=\sum_{i=1}^n\bar{J}_{ij},\; j=1,\ldots, n.$
\end{corollary}
Thus, in the protocols with additional uniform background links~\eqref{e_modeLK011115c}, consensus is determined by averaging the column entries of $(I+\tfrac1\delta L\big)^{-1},$ which is the parametric matrix of in-forests  of digraph~$\Gamma .$ Furthermore, if the arc weights of the added complete digraph tend to zero, then consensus 
is determined by the column means of the matrix of \emph{maximum\/} in-forests.

Comparing items\:2 of Theorem~\ref{t_depsop} and Corollary~\ref{c_limInf} we come to the following conclusion: the latent consensus of the protocol with uniform background links coincides with that of the protocol with hub uniformly subordinate to the agents.

\begin{theorem}\label{t_depsop}
$1.$ $\lim\limits_{\delta
\to+0}\lim\limits_{t\to\infty}\boldsymbol{y}(t) =
{\bf 1}'\big[\tilde{\boldsymbol{v}}^{\mathrm T}\!\bar{J}\;\; 0\big]
\boldsymbol{y}(0)$$;$

\noindent$2.$ If ${\boldsymbol{v}=\tfrac s n {\bf 1}},$ then
${\lim\limits_{\delta\to+0}\lim\limits_{t\to\infty}\boldsymbol{y}(t)
={\bf 1}'\big[\tfrac1n{\bf 1}^{\mathrm T}\!\bar{J}\;\:0\big]\boldsymbol{y}(0)
={\bf 1}'\,\tfrac1n\sum_{j=1}^n\bar{J}_{ \cdot j}y_j(0)},$
where ${\bar{J}_{\cdot j}=\sum_{i=1}^n\bar{J}_{ij}.}$
\end{theorem}

\section{Protocols of eigenprojection and orthogonal projection}
\label{s_ortho}

The same purpose, i.e., to develop a concept of latent consensus for the case where classical protocols give no consensus, can be achieved in another way. This can be done by means of the \emph{method of orthogonal projection}.
As distinct from the protocols of a dummy hub  and background links, where the dependency digraph $\Gamma $ is supplemented by additional links, in the method of orthogonal projection, the link structure is preserved, while the initial state is subject to a minimum necessary correction.
The resulting consensus has the form $\bar{J} S \boldsymbol{x}(0),$ where $S$ is the projection onto the \emph{consensus subspace\/} of the system determined by~$L.$
This consensus generally differs from those of protocol~\eqref{e_modeLK011115c} and the other protocols~\eqref{e_modeLK040715c} studied in this paper.

It should be noted that if the dependency digraph contains an in-tree, then all versions of latent consensus coincide with the ordinary
consensus~$\bar{J} \boldsymbol{x}(0).$

At last, let us consider the following second-order consensus protocol \cite{RenAtkins}:
\begin{equation}\label{150518a}
    \left[      \begin{array}{c}        \dot{\xi} \\        \dot{ \zeta} \\      \end{array}    \right]
    =\left[       \begin{array}{cc}
               0_{n\times n} & I_n \\
               -L & -\gamma L \\
             \end{array}           \right]  \left[          \begin{array}{c}            \xi \\             \zeta \\           \end{array}        \right],
    \end{equation}
where $\xi$ and $\zeta$ denote the position and velocity of the $i$th agent (vehicle).
An asymptotic consensus for this protocol is reached if $|\xi_i -\xi_j|\rightarrow 0$ and $|\zeta_i -\zeta_j|\rightarrow 0$ as $t\rightarrow \infty$.

We show that if $[\xi^{\rm T},\zeta^{\rm T}]^{\rm T} $ is a solution to the system of equations \eqref{150518a}, then for large $t$ we have
\begin{equation}\label{150518b}
    \left[      \begin{array}{c}        \dot{\xi}(t) \\        \dot{ \zeta}(t) \\      \end{array}    \right]
    =\left[             \begin{array}{cc}
             L^\vdash   & t L^\vdash   \\
                0_{n\times n} &   L^\vdash   \\
             \end{array}           \right].  \left[          \begin{array}{c}            \xi (0)\\             \zeta(0) \\           \end{array}        \right],
    \end{equation}
where $L^\vdash$ is the eigenprojection of $L$.

\section*{Acknowledgments}
This work was partially supported by the Presidium of RAS (Program No.~30 ``Theory and technologies of multi-level decentralized group control in the context of conflict and cooperation''). The work of PC was partially supported by the RSF (project no.\ 16-11-00063 granted to IRE RAS).

\bigskip\medskip
\noindent \textbf{Keywords:} consensus, multi-agent system, decentralized control,  DeGroot's iterative pooling, PageRank, Laplacian matrix of a digraph.

\bigskip
\noindent \textbf{AMS Subject Classification: 93A14, 68T42, 15B51, 05C50, 05C05, 60J22}

\end{document}